\documentclass[12pt,reqno, oneside]{amsart}
\usepackage{geometry}
\usepackage{xcolor}
\usepackage{ytableau}
\usepackage{hyperref}
\renewcommand{\l}{\lambda}
\newcommand{\inv}{\theta}

\definecolor{dred}{rgb}{.9,.1,0}
\definecolor{dgreen}{rgb}{0,.6,.1}
\definecolor{dblue}{rgb}{.1,.1,.8}

\begin{document}
\title[On the Schur function expansion]{On the  Schur function expansion\\
of a symmetric quasi-symmetric function}
\author{Ira M. Gessel$^*$}
\address{Department of Mathematics\\
   Brandeis University\\
   Waltham, MA 02453-2700}
\email{gessel@brandeis.edu}
\date{March 25, 2018}
\thanks{$^*$Supported by a grant from the Simons Foundation (\#427060, Ira Gessel)}
\begin{abstract}
Egge, Loehr, and Warrington proved a formula for the Schur function expansion of a symmetric function in terms of its expansion in fundamental quasi-symmetric functions. Their formula involves the coefficients of a modified inverse Kostka matrix. Recently Garsia and Remmel gave a simpler reformulation of
Egge, Loehr, and Warrington's result, with a new proof. We give here a simple proof of Garsia and Remmel's version, using a sign-reversing involution.

\end{abstract}
\maketitle

Egge, Loehr, and Warrington \cite{elw} proved a formula, involving the coefficients of a modified inverse Kostka matrix, for the Schur function expansion of a symmetric function in terms of its expansion as a linear combination of fundamental quasi-symmetric functions. We recall that for a composition $L=(L_1,\dots, L_k)$, the fundamental quasi-symmetric function $F_L$ is defined by
\begin{equation*}
F_L=\sum_{i_1,\dots, i_k} x_{i_1}\cdots x_{i_k},
\end{equation*}
where the sum is over all positive integers $i_1,\dots, i_k$ satisfying $i_1\le i_2\le\dots\le i_k$ and $i_j < i_{j+1}$ if $j\in\{L_1, L_1+L_2,\dots, L_1+L_2+\cdots +L_{k-1}\}$.

Garsia and Remmel \cite{gr} gave a simpler reformulation of Egge, Loehr, and Warrington's result.
For any composition $L$, we define the Schur function $s_L$ by the Jacobi-Trudi determinant of complete symmetric functions: $s_L = \det (h_{L_i-i+j})$, where 
$h_k$ is the complete symmetric function and $h_k=0$ for $k<0$. As explained below, for every composition $L$, $s_L$ is either an ordinary Schur function, the negative of an ordinary Schur function, or zero. 

Garsia and Remmel's reformulation is  that if $f$ is symmetric and $f=\sum_L c_L F_L$, then 
$f = \sum_L c_L s_L$. (There will usually be some cancellation in this formula.)

We give here a short combinatorial proof of Garsia and Remmel's reformulation. 
By linearity, it is sufficient to prove the formula for the case in which $f$ is a Schur function. We will show that for any partition $\l$, if $s_\l = \sum_L c_L F_L$ then $s_\l = \sum_L c_L s_L$.

For example, if $\l = (4,1)$ then $s_\l = F_{(4,1)}+F_{(3,2)}+F_{(2,3)}+F_{(1,4)}$.
We have $s_{(2,3)}=0$ and $s_{(1,4)}=-s_{(3,2)}$, so 
\[s_{(4,1)}+s_{(3,2)}+s_{(2,3)}+s_{(1,4)}=s_{(4,1)}+s_{(3,2)}+0-s_{(3,2)}=s_{(4,1)},\]
confirming the formula in this case.

Let $L=(L_1, \dots, L_k)$ be a composition. If $2\le i\le n$ and $L_{i}\ge2$, we define the composition $L^{(i)}$ to be  
$(L_1, \dots, L_{i-2}, L_{i}-1, L_{i-1} +1, L_{i+1}, \dots, L_k)$.
In other words, $L^{(i)}$ is obtained from $L$ by replacing $L_{i-1}, L_{i}$ with $L_{i}-1, L_{i-1}+1$, and leaving the other entries unchanged.
It follows from the Jacobi-Trudi determinant that
$ s_{L^{(i)}}=-s_L$. 
Thus if $L^{(i)}=L$ then $s_L = 0$ and an easy induction argument shows that for any composition $L$, $s_L$ is either 0 or $\pm s_\l$ for some partition $\l$. (This can also be shown be rearranging the rows of the Jacobi-Trudi determinant.)

The expansion of a Schur function into the fundamental quasi-symmetric functions is  a well-known consequence of Richard Stanley's theory of P-partitions \cite[Theorem 7.19.7, p.~361]{ec2}. A \emph{descent} of a standard tableau $T$ is an integer $i$ such that $i+1$ appears in a lower row in $T$ (in English notation) than $i$. Let the descents of the tableau $T$ with entries  $1, \ldots, n$ be $d_1< d_2<\dots< d_j$. The \emph{descent composition} of $T$, which we denote by $C(T)$, is the composition $(d_1, d_2 - d_1, \dots, d_j-d_{j-1}, n-d_j)$ of $n$. Then 
\begin{equation*}
s_\l = \sum_T F_{C(T)},
\end{equation*}
where the sum is over all standard tableaux $T$ of shape $\l$. So we need to prove that
\begin{equation}
\label{e-schursum}
s_\l = \sum_T s_{C(T)}.
\end{equation}

There is a unique standard tableau of shape $\l$ with descent composition $\l$, called the \emph{superstandard tableau}. It has entries $1, 2, \dots, \l_1$ in the first row, entries $\l_1+1, \l_1+2,\dots, \l_1+\l_2$ in the second row, and so on. If $T$ is superstandard of shape $\l$ then $C(T)=\l$, so $s_{C(T)}=s_\l$.

We will define a shape-preserving involution $\inv$ on standard but not superstandard tableaux,
with the property that $s_{C(\inv(T))}=-s_{C(T)}$. This property implies that if $\inv(T)=T$ then $s_{C(T)}=0$. So in the sum on the right side of \eqref{e-schursum} everything cancels except the term corresponding to the superstandard tableau of shape $\l$, which contributes $s_\l$, thus proving \eqref{e-schursum}.

Let $T$ be a standard tableau with descent set $S=\{d_1<d_2<\dots<d_j\}$. We define the \emph{$i$th run} of $T$, for $i$ from 1 to $j+1$, to be the skew subtableau of $T$ consisting of the elements $d_{i-1}+1, d_{i-1}+2,\dots,d_i$, where we set $d_0=0$ and $d_{j+1}=n$. Thus the number of elements in the $i$th run of $T$ is the $i$th part in the descent composition of~$T$.

For example, in the following tableau, the elements of the first run are colored red, the elements of the second run are colored green, and the elements of the third run are colored blue.
\[
\def\r{\color{dred}}
\def\g{\color{dgreen}}
\def\b{\color{dblue}}
\begin{ytableau}
\r1&\r2&\r3&\g6&\g7&\b9\\
\g4&\g5\\
\b8
\end{ytableau}
\]

We define the involution $\inv$ first for tableaux with exactly two runs.
If the tableau $T$ has two runs then the shape  of $T$ has two parts. There are  $\l_1-\l_2+1$ standard tableaux  of shape 
$(\l_1,\l_2)$ with two runs. Each such tableau is uniquely determined by an integer $j$ with $\l_2\le j \le \l_1$ for which
the first run contains $1, 2, \dots, j$, all in the first row, and the second run contains $j+1$, $j+2$, \dots, $\l_1+\l_2$, with $j+1,j+2,\dots, j+\l_2$ in the second row and $j+\l_2+1,\dots, \l_1+\l_2$ in the first row. 
Let $T_j$ be this tableau, where $\l = (\l_1, \l_2)$ is fixed. Then
the descent composition for $T_j$ is $(j, n-j)$, where $n=|\l|=\l_1+\l_2$. 

The superstandard tableau of shape $\l$ is $T_{\l_1}$. 
For $\l_2\le j \le \l_1-1$, we define $\theta(T_j)$ to be the tableau with descent composition $(n-j-1, j+1)$; i.e., $\theta(T_j) = T_{n-j-1}$.  In order for this definition to be valid, we must have $\l_2\le n-j-1\le \l_1-1$. For the first inequality we have $(n-j-1)-\l_2 = (\l_1+\l_2 -j-1)-\l_2 = (\l_1-1)-j\ge0$, and for the second inequality we have $(\l_1-1) -(n-j-1)=(\l_1-1)- (\l_1 +\l_2 -j-1) = j-\l_2\ge0$.
Thus $\theta$ is well-defined and  $s_{C(\inv(T_j))}=-s_{C(T_j)}$.

For example, if $\l=(4,2)$ and $j=3$ then $T_3$ is 
\[
\begin{ytableau}
1&2&3&6\\
4&5
\end{ytableau}
\]
with descent composition $(3,3)$, and $\theta(T_3)=T_2$ is
\[
\begin{ytableau}
1&2&5&6\\
3&4
\end{ytableau}
\]
with descent composition $(2,4)$.

Next, we define $\theta$ for tableaux $T$ of arbitrary shape $\l$ in which the first row is not $1, 2, \dots, \l_1$. Here we apply $\theta$ as defined above to the first two runs of $T$ and leave the rest of $T$ unchanged. So, for example, $\theta$ applied to 
\[
\begin{ytableau}
1&2&3&6&8&9\\
4&5&7
\end{ytableau}
\]
gives
\begin{equation*}
\begin{ytableau}
1&2&5&6&8&9\\
3&4&7
\end{ytableau}
\end{equation*}

\vspace{3pt}

Note that  we may extend $\theta$ as just defined in an obvious way to tableaux with any distinct entries, not necessarily $1, 2, \dots, n$.

In the general case suppose that the first $k$ rows of $T$ constitute a superstandard tableau but the first $k+1$ rows do not. (So $T$ must have at least $k+2$ rows.) Then to compute $\theta(T)$ we leave the first $k$ rows unchanged and apply  $\theta$ 
to the subtableau of $T$ consisting of rows $k+1$, $k+2$, \dots. It is clear that $\theta$ has the desired property: for every non-superstandard tableau $T$, we have $C(\theta(T)) = C(T)^{(i)}$ for some $i$, so 
$s_{C(\inv(T))}=-s_{C(T)}$.

For example, suppose that $T$ is the standard tableau
\begin{equation*}
\begin{ytableau}
1&2&3&4&5\\
6&7&9\\
8
\end{ytableau}
\end{equation*}
with descent composition $C(T)= (5,2,2)$. The first row of $T$ is superstandard but the first two rows are not. So $\theta(T)$ is 
\begin{equation*}
\begin{ytableau}
1&2&3&4&5\\
6&8&9\\
7
\end{ytableau}
\end{equation*}
with descent composition $(5,1,3) = C(T)^{(3)}$.


\begin{thebibliography}{9}
\bibitem{elw}
E. Egge, N. A. Loehr, and G. Warrington,
\emph{From quasisymmetric expansions to Schur expansions via a modified
   inverse Kostka matrix},
   European J. Combin.
   \textbf{31} (2010)
   2014--2027.
   
\bibitem{gr}
A. Garsia and J. Remmel, 
\emph{A note on passing from a quasi-symmetric function expansion to a Schur function expansion of a symmetric function}, \href{https://arxiv.org/abs/1802.09686}{\texttt{arXiv:1802.09686 [math.CO]}},  2018.

\bibitem{ec2}
R. P. Stanley,
\emph{Enumerative Combinatorics, Vol. 2},
Cambridge University Press, 1999.

\end{thebibliography}
\end{document}